\documentclass{amsart}

\usepackage{amssymb}
\usepackage{eepic}
\usepackage{color}

\allowdisplaybreaks

\numberwithin{equation}{section}

\newtheorem{theorem}{Theorem}[section]
\newtheorem{lemma}[theorem]{Lemma}

\newtheorem{remark}[theorem]{Remark}
\newtheorem{problem}[theorem]{Problem}

\newtheorem{TheoA}{Theorem A}

\newcommand{\Z}{\mathbb{Z}}
\newcommand{\R}{\mathbb{R}}
\newcommand{\C}{\mathbb{C}}

\newcommand{\summ}{\sum\nolimits}

\def\1{\mathbf{1}}
\def\H{\mathcal{H}}

\def\M{\mathcal{M}}

\newcommand{\dem}{\noindent {\bf Proof. }}

\newcommand{\demA}{\noindent {\bf Proof of Theorem A for $p < \infty$. }}
\newcommand{\demAA}{\noindent {\bf Proof of Theorem A for $p = \infty$. }}

\newcommand{\fin}{\hspace*{\fill} $\square$ \vskip0.2cm}

\def\mean{- \hskip-11pt \int}

\addtocounter{tocdepth}{-1}

\begin{document}

\title[Atomic blocks for martingales]{Atomic blocks for \\ noncommutative martingales}

\null

\vskip-15pt

\null

\author[Conde-Alonso and Parcet]
{Jose M. Conde-Alonso and Javier Parcet}

\maketitle

\begin{abstract}
Given a probability space $(\Omega,\Sigma,\mu)$, the Hardy space $\mathrm{H}_1(\Omega)$ which is associated to the martingale square function does not admit a classical atomic decomposition when the underlying filtration is not regular. In this paper we construct a decomposition of $\mathrm{H}_1(\Omega)$ into \lq atomic blocks\rq${}$ in the spirit of Tolsa, which we will introduce for martingales. We provide three proofs of this result. Only the first one also applies to noncommutative martingales, the main target of this paper. The other proofs emphasize alternative approaches for commutative martingales. One might be well-known to experts, using a weaker notion of atom and approximation by atomic filtrations. The last one adapts Tolsa's argument replacing medians by conditional medians.
\end{abstract}

\addtolength{\parskip}{+1ex}

\section*{{\bf Introduction}}

Let $(\Omega, \Sigma,\mu)$ be a probability space equipped with a filtration $(\Sigma_k)_{k \ge 1}$ whose union generates $\Sigma$. Let us write $\mathsf{E}_k$ to denote the conditional expectation onto $\Sigma_k$-measurable functions and $\Delta_k = \mathsf{E}_k - \mathsf{E}_{k-1}$ for the associated differences, with the convention that $\Delta_1 = \mathsf{E}_1$. Given $f \in L_1(\Omega)$, we shall usually write $f_k$ and $df_k$ for $\mathsf{E}_k f$ and $\Delta_k f$ respectively. Once the filtration $(\Sigma_k)_{k \ge 1}$ is fixed, the martingale Hardy space $\mathrm{H}_1(\Omega)$ is the subspace of functions $f$ in $L_1(\Omega)$ whose $\mathrm{H}_1(\Omega)$-norm defined below is finite $$\|f\|_{\mathrm{H}_1(\Omega)} \, = \, \Big\| \Big( \sum_{k \ge 1} |df_k|^2\Big)^\frac12 \Big\|_1.$$ As it was proved by Davis \cite{D}, we obtain an equivalent norm after replacing the martingale square function above by Doob's martingale maximal function. On the contrary, replacing the martingale square function by its conditioned form we get the so-called little Hardy space $\mathrm{h}_1(\Omega)$. In other words, the subspace of functions $f$ in $L_1(\Omega)$ whose $\mathrm{h}_1(\Omega)$-norm below is finite under the convention $\mathsf{E}_{k-1} |df_k|^2 = |f_1|^2$ when $k=1$ $$\|f\|_{\mathrm{h}_1(\Omega)} \, = \, \Big\| \Big( \sum_{k \ge 1} \mathsf{E}_{k-1} |df_k|^2 \Big)^\frac12 \Big\|_1.$$ Both spaces are fair generalizations of the Euclidean Hardy space. Namely, if we pick the standard dyadic filtration in $\R^n$, it turns out that $\mathrm{H}_1(\Omega)$ is by all means the dyadic form of $\mathrm{H}_1$, whereas we have $\mathrm{h}_1(\Omega) \simeq \mathrm{H}_1(\Omega)$ for regular filtrations as it happens in the dyadic setting. It is in the case of nonregular filtrations when both spaces have their own identity. In general, we have $\mathrm{h}_1(\Omega) \subsetneq \mathrm{H}_1(\Omega)$ and more precisely we know from \cite{JX} that $$\|f\|_{\mathrm{H}_1(\Omega)} \, \sim \, \inf_{f = g + h} \, \|g\|_{\mathrm{h}_1(\Omega)} + \sum_{k \ge 1} \|dh_k\|_1.$$ We refer to Garsia's book \cite{G} for more information on martingale Hardy spaces. 

Given $1 < p \le \infty$, a function $\mathrm{a}: \Omega \to \C$ is called a martingale $p$-atom when $\mathrm{a}$ is $\Sigma_1$-measurable and $\|\mathrm{a}\|_1 = 1$ or there exists $k \ge 1$ and a measurable set $\mathrm{A} \in \Sigma_k$ such that 
\begin{itemize}
\item $\mathsf{E}_k(\mathrm{a}) = 0$,

\vskip2pt

\item $\mathrm{supp} (\mathrm{a}) \subset \mathrm{A}$, 

\item $\|\mathrm{a}\|_p \le \mu(\mathrm{A})^{- \frac{1}{p'}}$ for $\frac1p + \frac{1}{p'}=1$.
\end{itemize}
The motivation for this article is the fact that (with this notion of atom) no atomic description is known for the space $\mathrm{H}_1(\Omega)$. On the contrary, $\mathrm{h}_1(\Omega)$ always admits an atomic decomposition. Indeed, define the atomic Hardy spaces as
\begin{eqnarray*}
\mathrm{h}_{\mathrm{at}}^1(\Omega) & = & \Big\{ f \in L_1(\Omega) \ \big| \ f = \sum_{j \ge 1} \lambda_j \mathrm{a}_j, \ \mathrm{a}_j \ \mbox{$2$-atom}, \ \sum_{j \ge 1} |\lambda_j| < \infty \Big\}, \\ \mathrm{h}_{\mathrm{at},p}^1(\Omega) & = & \Big\{ f \in L_1(\Omega) \ \big| \ f = \sum_{j \ge 1} \lambda_j \mathrm{a}_j, \ \mathrm{a}_j \ \mbox{$p$-atom}, \ \sum_{j \ge 1} |\lambda_j| < \infty \Big\}.
\end{eqnarray*}
The norm is the infimum of $\sum_j |\lambda_j|$ over all decompositions of $f$. As a combination of \cite{H,W}, we know that $\mathrm{h}_1(\Omega) \simeq \mathrm{h}_{\mathrm{at},p}^1(\Omega)$ for $1 < p \le \infty$. This yields an atomic decomposition of $\mathrm{h}_1(\Omega)$, which works for $\mathrm{H}_1(\Omega)$ when the filtration is regular.  

Atomic decompositions are useful to provide endpoint estimates for singular operators $T$ failing to be bounded in $L_1(\Omega)$. Indeed, this typically reduces ---under mild regularity assumptions--- to bound uniformly the $L_1$-norm of $T(\mathrm{a})$ for an arbitrary atom $\mathrm{a}$, which is easier than proving the $\mathrm{H}_1 \to L_1$ boundedness of $T$ due to the particular structure of atoms. The drawback of the martingale atoms described above is that they are useless for $\mathrm{H}_1(\Omega)$ when the filtration is not regular. This is significant because in that case the spaces $\mathrm{h}_1(\Omega)$ are not endpoint interpolation spaces in the $L_p$ scale, whereas the spaces $\mathrm{H}_1(\Omega)$ are. Therefore, the goal of this paper is to provide an alternative atomic decomposition for $\mathrm{H}_1(\Omega)$ suitable for arbitrary filtrations, and also for classical and noncommutative martingales. 

Our approach is strongly motivated by the work of Tolsa on the so-called RBMO spaces \cite{To}. Namely, it is well-known that we have $\mathrm{h}_1(\Omega)^* \simeq \mathrm{bmo}(\Omega)$ and also $\mathrm{H}_1(\Omega)^* \simeq \mathrm{BMO}(\Omega)$ where both martingale BMO spaces are respectively defined as the functions $f$ in $L_2(\Omega)$ with finite norm
\begin{eqnarray*}
\|f\|_{\mathrm{bmo}(\Omega)} & = & \sup_{k \ge 1} \Big\| \Big( \mathsf{E}_k \big| f - \mathsf{E}_{k} f \big|^2 \Big)^\frac12 \Big\|_\infty, \\ \|f\|_{\mathrm{BMO}(\Omega)} & = & \sup_{k \ge 1} \Big\| \Big( \mathsf{E}_k \big| f - \mathsf{E}_{k-1} f \big|^2 \Big)^\frac12 \Big\|_\infty. 
\end{eqnarray*}
It is easily checked that we have the norm equivalence $$\|f\|_{\mathrm{BMO}(\Omega)} \, \simeq \, \|f\|_{\mathrm{bmo}(\Omega)} + \sup_{k \ge 1} \|df_k\|_\infty.$$ In analogy, Tolsa's RBMO norm is the sum of a \lq doubling\rq${}$ BMO norm plus a term which measures the \lq distance\rq${}$ between averages over nested pairs of doubling cubes. This viewpoint is fruitful in both directions. Indeed, nondoubling techniques are adapted here (in one of the approaches we follow) for martingales whereas martingale techniques are used in \cite{CMP} for nondoubling spaces.

Tolsa's construction of the predual of RBMO is therefore our model to produce an atomic type decomposition of $\mathrm{H}_1(\Omega)$. A $\Sigma$-measurable function $\mathrm{b}: \Omega \to \C$ will be called a \emph{martingale $p$-atomic block} when $\mathrm{b} \in L_1(\Omega,\Sigma_1,\mu)$ or there exists $k \ge 1$ such that the following properties hold  
\begin{itemize}
\item $\mathsf{E}_k(\mathrm{b}) = 0$,

\vskip2pt

\item $\mathrm{b} = \sum_j \lambda_j \mathrm{a}_j$ where 

\vskip4pt

\begin{itemize}
\item[$\circ$] $\mathrm{supp}(\mathrm{a}_j) \subset \mathrm{A}_j$, 

\item[$\circ$] $\|\mathrm{a}_j\|_p \le \mu(\mathrm{A}_j)^{-\frac{1}{p'}} \frac{1}{k_j - k + 1}$,
\end{itemize}

\vskip3pt

\noindent for certain $k_j \ge k$ and $\mathrm{A}_j \in \Sigma_{k_j}$. Call each such $\mathrm{a}_j$ a \emph{$p$-subatom}.
\end{itemize}
Given a martingale $p$-atomic block, set $$|\mathrm{b}|_{\mathrm{atb},p}^1 \, = \, \begin{cases} \displaystyle \int_\Omega |\mathrm{b}(\omega)| \, d\mu(\omega) & \mbox{when } \ \mathrm{b} \in L_1(\Omega, \Sigma_1, \mu), \\ \displaystyle  \inf_{\begin{subarray}{c} \mathrm{b} = \sum_j \lambda_j \mathrm{a}_j \\ \mathrm{a}_j \ p-\mathrm{subatom} \end{subarray}} \ \sum_{j \ge 1} |\lambda_j| & \mbox{when } \ \mathrm{b} \notin L_1(\Omega, \Sigma_1, \mu). \end{cases}$$ Then we define the atomic block Hardy spaces 
\begin{eqnarray*}
\mathrm{H}_{\mathrm{atb}}^1(\Omega) & = & \Big\{ f \in L_1(\Omega) \ \big| \ f = \summ_i \mathrm{b}_i, \ \mathrm{b}_i \ \mbox{martingale } \mbox{$2$-atomic block} \Big\}, \\ \mathrm{H}_{\mathrm{atb},p}^1(\Omega) & = & \Big\{ f \in L_1(\Omega) \ \big| \ f = \summ_i \mathrm{b}_i, \ \mathrm{b}_i \ \mbox{martingale } \mbox{$p$-atomic block} \Big\},
\end{eqnarray*}
which come equipped with the norm $$\|f\|_{\mathrm{H}_{\mathrm{atb},p}^1(\Omega)} \, = \, \inf_{\begin{subarray}{c} f = \sum_i \mathrm{b}_i \\ \mathrm{b}_i \ p-\mathrm{atomic \, block} \end{subarray}} \sum_{i \ge 1} |\mathrm{b}_i|_{\mathrm{atb},p}^1 \, = \, \inf_{\begin{subarray}{c} f = \sum_i \mathrm{b}_i \\ \mathrm{b}_i = \sum_j \lambda_{ij} \mathrm{a}_{ij} \end{subarray}} \ \sum_{i,j \ge 1} |\lambda_{ij}|,$$ where the $\mathrm{a}_{ij}$'s above are taken to be $p$-subatoms of $\mathrm{b}_i$. Note that $\lambda_{ij} = \delta_{j1} \|\mathrm{b}_i\|_1$ for atomic blocks $\mathrm{b}_i \in L_1(\Omega,\Sigma_1,\mu)$. With this definition of atomic blocks, $\mathrm{H}_1 \to \mathrm{L}_1$ boundedness reduces to $$\|T(\mathrm{b})\|_1 \, \le \, \mathrm{c}_0 |\mathrm{b}|_{\mathrm{atb},p}^1$$ under mild regularity conditions for some $\mathrm{c}_0$ independent of the $p$-atomic block $\mathrm{b}$.

\begin{TheoA}
There exists an isomorphism $$\mathrm{H}_1(\Omega) \, \simeq \, \mathrm{H}_{\mathrm{atb},p}^1(\Omega) \quad \mbox{for} \quad 1 < p \le \infty.$$ In fact, an analogous result holds also for noncommutative martingales. 
\end{TheoA}

We have deliberately omitted the definition of atomic block for noncommutative martingales, which is postponed to Section \ref{Sect2}. We shall only provide one proof of Theorem A which is valid for noncommutative martingales, although two additional arguments will be given in the commutative setting. Our main proof is perhaps the simplest one, relying on a noncommutative form of Davis decomposition from \cite{JMe,Pe}. An alternative proof exploits a weaker notion of atom which might be folklore or at least well-known to experts. It however requires to approximate general filtrations by atomic ones, something which seems to be out of the scope in the noncommutative setting. Our last proof avoids such approximation argument adapting Tolsa's argument \cite{To} with conditional medians instead of medians. Our noncommutative results are in line with \cite{BCP,HM,Pe}.

\section{{\bf Noncommutative martingales}} \label{Sect2}

The theory of noncommutative martingale inequalities started with Cuculescu \cite{C}, but it did not receive significant attention until the work of Pisier/Xu \cite{PX} about the noncommutative analogue of Burkholder/Gundy inequalities. After it, most of the classical results on martingale $L_p$ inequalities have found a noncommutative analogue, see \cite{HM,J,JMe,JMu,JPe,JX,Mu,PR,Pe,X} and the references therein for basic definitions and results. Here we shall just introduce martingale $p$-atomic blocks and related notions in the noncommutative setting.      

A noncommutative probability space is a pair $(\M,\tau)$ 
formed by a von Neumann algebra $\M$ and a normal faithful finite trace $\tau$, normalized so that $\tau(\1_\M) = 1$ for the unit $\1_\M$ of $\M$. A filtration in $\M$ is an increasing sequence $(\M_k)_{k \ge 1}$ of von Neumann subalgebras of $\M$ satisfying that their union is weak-$*$ dense in $\M$. Assume there exists a normal conditional expectation $$\mathsf{E}_k: \M \to \M_k$$ for every $k \ge 1$. Each $\mathsf{E}_k$ is trace preserving, unital and completely positive. In particular, $\mathsf{E}_k: L_p(\M) \to L_p(\M_k)$ defines a contraction for $1 \le p \le \infty$. These maps satisfy the bimodule property $\mathsf{E}_k(\alpha f \beta) = \alpha \mathsf{E}_k(f) \beta$ for $\alpha, \beta \in \M_k$. If we set $\Delta_k = \mathsf{E}_k - \mathsf{E}_{k-1}$ and write $\mathsf{E}_kf = f_k$ and $\Delta_k f = df_k$ for $f \in L_1(\M)$ (as in the commutative setting) then $\mathrm{H}_1(\M)$ is defined as the subspace of operators $f \in L_1(\M)$ with finite norm $$\|f\|_{\mathrm{H}_1(\M)} \, = \, \inf_{\begin{subarray}{c} f = g+h \\ g,h \in L_1(\M) \end{subarray}} \|g\|_{\mathrm{H}_1^\mathrm{c}(\M)} + \|h^*\|_{\mathrm{H}_1^\mathrm{c}(\M)},$$ where the column Hardy norm is given by $$\|f\|_{\mathrm{H}_1^\mathrm{c}(\M)} \, = \, \Big\| \Big( \sum_{k \ge 1} df_k^* df_k \Big)^\frac12 \Big\|_1.$$ The little Hardy space is defined similarly with $$\|f\|_{\mathrm{h}_1^\mathrm{c}(\M)} \, = \, \Big\| \Big( \sum_{k \ge 1} \mathsf{E}_{k-1} \big( df_k^* df_k \big) \Big)^\frac12 \Big\|_1.$$ On the other hand, $\mathrm{BMO}(\M)$ is the subspace of $L_2(\M)$ with $$\|f\|_{\mathrm{BMO}(\M)} \, = \, \max \Big\{ \|f\|_{\mathrm{BMO}_\mathrm{c}(\M)}, \|f^*\|_{\mathrm{BMO}_\mathrm{c}(\M)} \Big\}$$ where the column BMO norm is given by the following expression $$\|f\|_{\mathrm{BMO}_\mathrm{c}(\M)} \, = \, \sup_{k \ge 1} \Big\| \Big(\mathsf{E}_k \big( (f - \mathsf{E}_{k-1}f)^* (f - \mathsf{E}_{k-1}f) \big) \Big)^\frac12 \Big\|_\M.$$ Of course, $\mathrm{bmo}(\M)$ arises when we replace $\mathsf{E}_{k-1}$ by $\mathsf{E}_k$ in the identity above. 

We are now ready to define martingale $p$-atomic blocks in the noncommutative setting. As expected, we find row and column forms of these objects. We will say that an (unbounded) operator $\mathrm{b}$ affiliated with the von Neumann algebra $\M$ is a \emph{column martingale $p$-atomic block} when $\mathrm{b} \in L_1(\M_1,\tau)$ or there exists an index $k \ge 1$ such that
\begin{itemize}
\item $\mathsf{E}_k(\mathrm{b}) = 0$,

\vskip2pt

\item $\mathrm{b} = \sum_j \lambda_j \mathrm{a}_j$ where

\vskip4pt

\begin{itemize}
\item[$\circ$] $\mathrm{a}_j q_j = \mathrm{a}_j$, 

\item[$\circ$] $\|\mathrm{a}_j\|_p \le \tau(q_j)^{-\frac{1}{p'}} \frac{1}{k_j - k + 1}$,
\end{itemize}

\vskip3pt

\noindent for some $k_j \ge k$ and projections $q_j \in \M_{k_j}$. 
\end{itemize}
Each such $\mathrm{a}_j$ will be called a \emph{column $p$-subatom}. Similarly, \emph{row $p$-atomic blocks} are defined when the support identity $q_j \mathrm{a}_j = \mathrm{a}_j$ holds instead. In particular, both conditions hold for self-adjoint atomic blocks. Given a column $p$-atomic block $\mathrm{b}$ set $$|\mathrm{b}|_{\mathrm{atb},p}^{1,\mathrm{c}} = \tau(|\mathrm{b}|)$$ when $\mathrm{b} \in L_1(\M_1,\tau)$ and otherwise $$|\mathrm{b}|_{\mathrm{atb},p}^{1,\mathrm{c}} \, = \, \inf_{\begin{subarray}{c} \mathrm{b} = \sum_j \lambda_j \mathrm{a}_j \\ \mathrm{a}_j \ p-\mathrm{subatom} \end{subarray}} \ \sum_{j \ge 1} |\lambda_j|.$$ Then we define the atomic block Hardy spaces 
\begin{eqnarray*}
\mathrm{H}_{\mathrm{atb}}^{1,\mathrm{c}}(\M) & = & \Big\{ f \in L_1(\M) \ \big| \ f = \summ_i \mathrm{b}_i, \ \mathrm{b}_i \ \mbox{column } \mbox{$2$-atomic block} \Big\}, \\ \mathrm{H}_{\mathrm{atb},p}^{1,\mathrm{c}}(\M) & = & \Big\{ f \in L_1(\M) \ \big| \ f = \summ_i \mathrm{b}_i, \ \mathrm{b}_i \ \mbox{column } \mbox{$p$-atomic block} \Big\},
\end{eqnarray*}
which come equipped with the following norm $$\|f\|_{\mathrm{H}_{\mathrm{atb},p}^{1,\mathrm{c}}(\M)} \, = \, \inf_{\begin{subarray}{c} f = \sum_i \mathrm{b}_i \\ \mathrm{b}_i \ p-\mathrm{atomic \, block} \end{subarray}} \sum_{i \ge 1} |\mathrm{b}_i|_{\mathrm{atb},p}^{1,\mathrm{c}} \, = \, \inf_{\begin{subarray}{c} f = \sum_i \mathrm{b}_i \\ \mathrm{b}_i = \sum_j \lambda_{ij} \mathrm{a}_{ij} \end{subarray}} \ \sum_{i,j \ge 1} |\lambda_{ij}|,$$ where the $\mathrm{a}_{ij}$'s above are taken to be $p$-subatoms of $\mathrm{b}_i$. As in the commutative case, we pick $\lambda_{ij} = \delta_{j1} \|\mathrm{b}_i\|_1$ for atomic blocks $\mathrm{b}_i \in L_1(\M_1,\tau)$. Before stating the analogue of Theorem A for noncommutative martingales, we shall need the following  approximation lemma to legitimate our duality argument below.

\begin{lemma} \label{DensityLemma}
Given $\varepsilon > 0$ and $$f \in \mathrm{H}_{\mathrm{atb},p}^{1,\mathrm{c}}(\M),$$ there exist a finite family $(\mathrm{b}_i(\varepsilon))_{i \le \mathrm{M}}$ of column $p$-atomic blocks with 
\begin{itemize}
\item[i)] $\mathrm{b}_i(\varepsilon) \in L_p(\M)$,

\vskip7pt

\item[ii)] $\displaystyle \Big\| f - \sum_{i=1}^{\mathrm{M}} \mathrm{b}_i(\varepsilon) \Big\|_{\mathrm{H}_{\mathrm{atb},p}^{1,\mathrm{c}}(\M)} < \varepsilon$.

\item[iii)] $\displaystyle \sum_{i=1}^{\mathrm{M}} |\mathrm{b}_i(\varepsilon)|_{\mathrm{atb},p}^{1,\mathrm{c}} < \Big\| \sum_{i=1}^{\mathrm{M}} \mathrm{b}_i(\varepsilon) \Big\|_{\mathrm{H}_{\mathrm{atb},p}^{1,\mathrm{c}}(\M)} + \varepsilon$.
\end{itemize}
\end{lemma}

\dem Let $f = \sum_i \mathrm{b}_i$ be such that 
\begin{eqnarray*}
\Big\| f - \sum_{i=1}^{\mathrm{M}} \mathrm{b}_i \Big\|_{\mathrm{H}_{\mathrm{atb},p}^{1,\mathrm{c}}(\M)} & < & \delta, \\ \Big| \|f\|_{\mathrm{H}_{\mathrm{atb},p}^{1,\mathrm{c}}(\M)} - \sum_{i = 1}^{\mathrm{M}} |\mathrm{b}_i|_{\mathrm{atb},p}^{1,\mathrm{c}} \Big| & < & \delta,
\end{eqnarray*}
with $\delta = \delta(\varepsilon)$ small and $\mathrm{M} = \mathrm{M}(\delta)$ large enough. From these properties it is clear that all the assertions in the statement will follow as long as we can show that every column $p$-atomic block $\mathrm{b}$ can be $\delta$-approximated by another column $p$-atomic block $\mathrm{b}'$ living in $L_p(\M)$. Indeed, when $\mathrm{b} \in L_1(\M_1,\tau)$ it suffices to select an element $\mathrm{b}' \in L_p(\M_1,\tau) \subset L_1(\M_1,\tau)$ with $\|\mathrm{b} - \mathrm{b}'\|_1 < \delta$. Otherwise $$\mathrm{b} = \summ_j \lambda_j \mathrm{a}_j \quad \mbox{with} \quad \mathsf{E}_k(\mathrm{b})=0$$ is a sum of column $p$-subatoms. In that case, set $\mathrm{N} = \mathrm{N}(\delta)$ so that $$\sum_{j > \mathrm{N}} |\lambda_j| < \frac{\delta}{2k}$$ and define $$\mathrm{b}' = \sum_{j \le \mathrm{N}} \lambda_j \mathrm{a}_j + \mathsf{E}_1 \Big( \sum_{j > \mathrm{N}} \lambda_j \mathrm{a}_j \Big) =: \sum_{j \le \mathrm{N}+1} \lambda_j' \mathrm{a}_j'.$$  According to the definition of column $p$-atomic block, the following holds 
\begin{itemize}
\item $\mathsf{E}_1(\mathrm{b}') = \mathsf{E}_1(\mathrm{b}) = \mathsf{E}_1 \mathsf{E}_k (\mathrm{b}) = 0$,

\item If $(k_j',q_j') = (k_j,q_j)$ for $j \le \mathrm{N}$ and $(k_{\mathrm{N}+1}',q_{\mathrm{N}+1}') = (1,\1_\M)$ $$\mathrm{a}_j' q_j' = \mathrm{a}_j' \quad \mbox{and} \quad \|\mathrm{a}_j'\|_p \le k \, \tau(q_j')^{-\frac{1}{p'}} \mbox{$\frac{1}{k_j-1+1}$}$$ provided we normalize $\mathrm{a}_{\mathrm{N}+1}'$ so that $\lambda_{\mathrm{N}+1}' = \| \mathsf{E}_1 (\sum_{j > \mathrm{N}} \lambda_j \mathrm{a}_j) \|_1$.
\end{itemize}
This shows that $\mathrm{b}'$ is a column $p$-atomic block. Moreover $$\|\mathrm{b}'\|_p = \Big\| \sum_{j \le \mathrm{N}} \lambda_j \mathrm{a}_j - \mathsf{E}_1 \Big( \sum_{j \le \mathrm{N}} \lambda_j \mathrm{a}_j \Big) \Big\|_p \le 2 \sum_{j \le \mathrm{N}} |\lambda_j| \|\mathrm{a}_j\|_p < \infty.$$ Therefore, it just remains to prove the following estimate $$\big\| \mathrm{b} - \mathrm{b}'  \big\|_{\mathrm{H}_{\mathrm{atb},p}^{1,\mathrm{c}}(\M)} < \delta.$$ To that aim we identify $\mathrm{b} - \mathrm{b}'$ as a column $p$-atomic block $$\mathrm{b} - \mathrm{b}' = \sum_{j > \mathrm{N}} \lambda_j \mathrm{a}_j - \mathsf{E}_1 \Big( \sum_{j > \mathrm{N}} \lambda_j \mathrm{a}_j \Big) =: \sum_{j > \mathrm{N}} \widetilde{\lambda}_j \widetilde{\mathrm{a}}_j + \widetilde{\lambda}_{\mathrm{N}} \widetilde{\mathrm{a}}_{\mathrm{N}}$$ with $\widetilde{\mathrm{a}}_\mathrm{N}$ normalized so that $\widetilde{\lambda}_{\mathrm{N}} = \| \mathsf{E}_1 (\sum_{j > \mathrm{N}} \lambda_j \mathrm{a}_j) \|_1$. Then we find 
\begin{itemize}
\item $\mathsf{E}_1(\mathrm{b} - \mathrm{b}') = 0$,

\item If $(\widetilde{k}_j,\widetilde{q}_j) = (k_j,q_j)$ for $j > \mathrm{N}$ and $(\widetilde{k}_{\mathrm{N}},\widetilde{q}_{\mathrm{N}}) = (1,\1_\M)$ $$\widetilde{\mathrm{a}}_j \widetilde{q}_j = \widetilde{\mathrm{a}}_j \quad \mbox{and} \quad \|\widetilde{\mathrm{a}}_j\|_p \le k \, \tau(\widetilde{q}_j)^{-\frac{1}{p'}} \mbox{$\frac{1}{k_j-1+1}$}.$$ 
\end{itemize}
This makes it quite simple to estimate the $\mathrm{H}_{\mathrm{atb},p}^{1,\mathrm{c}}(\M)$-norm of $\mathrm{b} - \mathrm{b}'$ 
\begin{eqnarray*}
\big\| \mathrm{b} - \mathrm{b}' \big\|_{\mathrm{H}_{\mathrm{atb},p}^{1,\mathrm{c}}(\M)} & \le & \big| \mathrm{b} - \mathrm{b}' \big|_{\mathrm{atb},p}^{1,\mathrm{c}} \\ [4pt] & \le & k \Big[ \sum_{j > \mathrm{N}}  |\lambda_j| + \Big\| \mathsf{E}_1\Big( \sum_{j > \mathrm{N}}  \lambda_j \mathrm{a}_j \Big) \Big\|_1 \Big] \\ & \le & k \Big[ \sum_{j > \mathrm{N}}  |\lambda_j| + \sum_{j > \mathrm{N}}  |\lambda_j| \| \mathrm{a}_j \|_1 \Big] \ \le \ 2k \sum_{j > \mathrm{N}} |\lambda_j| \ < \ \delta.
\end{eqnarray*}
Here we used the inequality $\|\mathrm{a}_j\|_1 = \|\mathrm{a}_j q_j\|_1 \le \|\mathrm{a}_j\|_p \tau(q_j)^{\frac{1}{p'}} \le \frac{1}{k_j-k+1} \le 1$. \fin

\begin{theorem}
There exists an isomorphism $$\mathrm{H}_1^{\mathrm{c}}(\M) \, \simeq \, \mathrm{H}_{\mathrm{atb},p}^{1,\mathrm{c}}(\M) \quad \mbox{for} \quad 1 < p \le \infty.$$ In particular, we find the atomic block decomposition $\mathrm{H}_1(\M) \simeq \mathrm{H}_{\mathrm{atb},p}^1(\M)$.
\end{theorem}

\dem We need to show 
\begin{itemize}
\item[i)] $\mathrm{H}_{\mathrm{atb},p}^{1,\mathrm{c}}(\M) \subset \mathrm{H}_1^{\mathrm{c}}(\M)$,

\item[ii)] $\mathrm{H}_1^{\mathrm{c}}(\M) \subset \mathrm{H}_{\mathrm{atb},p}^{1,\mathrm{c}}(\M)$.
\end{itemize}
\textbf{Step 1.} For the first continuous inclusion we shall prove $$\mathrm{BMO}_\mathrm{c}(\M) \subset \mathrm{H}_{\mathrm{atb},p}^{1,\mathrm{c}}(\mathcal{M})^*,$$ which suffices by duality. Assume that $\phi \in \mathrm{BMO}_\mathrm{c}(\mathcal{M})$. Since $\phi \in L_{p'}(\mathcal{M})$ for any $1 < p \le \infty$, we may represent $\phi$ as a linear functional $L_\phi$ on $L_p(\mathcal{M})$ by the formula $$L_\phi(f) = \tau(f\phi^*).$$ According to Lemma \ref{DensityLemma}, it suffices to show that $$|L_\phi(f)| \, \le \, \|f\|_{\mathrm{H}_{\mathrm{atb},p}^{1,\mathrm{c}}(\M)} \|\phi\|_{\mathrm{BMO}_\mathrm{c}(\M)}$$ for every $f$ which can be written as a finite sum $f = \sum_i \mathrm{b}_i$ of column $p$-atomic blocks $\mathrm{b}_i \in L_p(\M)$. This clearly allows us to estimate $|L_\phi(f)| \le \sum_i |L_\phi(\mathrm{b}_i)|$ with the right hand side well-defined. In particular, it is enough to show that $$|L_\phi(\mathrm{b})| \lesssim |\mathrm{b}|_{\mathrm{atb},p}^{1,\mathrm{c}} \|\phi\|_{ \mathrm{BMO}_\mathrm{c}(\mathcal{M})}$$ for column $p$-atomic blocks $\mathrm{b} \in L_p(\M)$. When $\mathrm{b} \in L_p(\M_1,\tau)$ $$| L_\phi(\mathrm{b}) | \le \|\mathrm{b}\|_1 \|\mathsf{E}_1 \phi\|_{\infty} \le |\mathrm{b}|_{\mathrm{atb},p}^{1,\mathrm{c}} \|\phi\|_{ \mathrm{BMO}_\mathrm{c}(\mathcal{M})}.$$ Otherwise, we write $\mathrm{b} = \sum_j \lambda_j \mathrm{a}_j$ with $\mathsf{E}_k(\mathrm{b})=0$ and such that $$\mathrm{a}_j q_j = \mathrm{a}_j \quad , \quad \|\mathrm{a}_j\|_p \le \tau(q_j)^{-\frac{1}{p'}} \mbox{$\frac{1}{k_j-k+1}$}$$ for some $k_j \ge k$ and some projection $q_j \in \M_{k_j}$. Then we find that $$|L_\phi(\mathrm{b})| = |\tau(\mathrm{b} \phi^*) | = \big| \tau \big( \mathrm{b} (\phi - \mathsf{E}_k\phi)^* \big) \big| \le \summ_j |\lambda_j| \big\| \mathrm{a}_j(\phi - \mathsf{E}_k\phi)^* \big\|_1 =: \summ_j |\lambda_j| \mathrm{A}_j.$$ Hence, it remains to prove that $\sup_j \mathrm{A}_j \lesssim \|\phi\|_{\mathrm{BMO}_\mathrm{c}(\M)}$, which follows from
\begin{eqnarray*}
\mathrm{A}_j & \le & \|a_j\|_p \big\| (\phi - \mathsf{E}_k \phi) q_j \big\|_{p'} \\ [8pt] & \le & \tau(q_j)^{-\frac{1}{p'}} \frac{1}{k_j-k+1} \big\| (\phi - \mathsf{E}_k \phi) q_j \big\|_{p'} \\
& \le & \tau(q_j)^{-\frac{1}{p'}} \big\| (\phi - \mathsf{E}_{k_j} \phi) q_j \big\|_{p'} + \frac{1}{k_j-k+1}  \sum_{s=k+1}^{k_j} \|d\phi_s\|_\infty \ = \ \mathrm{B}_j + \mathrm{C}_j
\end{eqnarray*}
Indeed, this yields the estimate $$\mathrm{B}_j + \mathrm{C}_j \le \|\phi\|_{\mathrm{bmo}_\mathrm{c}(\M)} + \sup_{k \ge 1} \|d\phi_k\|_\infty \, \simeq \, \|\phi\|_{\mathrm{BMO}_\mathrm{c}(\M)}$$ where the inequality $\mathrm{B}_j \le \|\phi\|_{\mathrm{bmo}_\mathrm{c}(\M)}$ follows from Hong/Mei formulation of the John-Nirenberg inequality for noncommutative martingales \cite{HM}. In particular, this completes the proof of Step 1. 

\vskip5pt

\noindent \textbf{Step 2.} We now prove the inclusion $$\mathrm{H}_1^{\mathrm{c}}(\M) \subset \mathrm{H}_{\mathrm{atb},p}^{1,\mathrm{c}}(\M)$$ directly, without using duality. Here we would like to thank Marius Junge for suggesting us the noncommutative Davis decomposition (used below) as a possible tool in proving this inclusion. Let $f \in \mathrm{H}_1^\mathrm{c}(\M)$, by the noncommutative form of Davis decomposition \cite{Pe} we know that $f$ can be decomposed as $f = f_\mathrm{c} + f_\mathrm{d}$, where $$(f_\mathrm{c}, f_\mathrm{d}) \in \mathrm{h}_{\mathrm{at,c}}^1(\M) \times \mathrm{h}_{\mathrm{diag}}^1(\mathcal{M}).$$ On the other hand, since a column $p$-atom in the sense of \cite{BCP,HM} is in particular a column $p$-atomic block in our sense, we immediately find the following inequality $$\|f_\mathrm{c}\|_{\mathrm{H}_{\mathrm{atb},p}^{1,\mathrm{c}}(\mathcal{M})} \, \lesssim \, \|f\|_{\mathrm{H}_1^\mathrm{c}(\M)}.$$ The diagonal norm of $f_\mathrm{d}$ is given by $$\|f_\mathrm{d}\|_{\mathrm{h}_{\mathrm{diag}}^1(\M)} \, = \, \sum_{k \ge 1} \big\| \Delta_k(f_\mathrm{d}) \big\|_1 \, \lesssim \, \|f\|_{\mathrm{H}_1^\mathrm{c}(\M)}.$$ Therefore, the goal is to show that we have $$\|f_\mathrm{d}\|_{\mathrm{H}_{\mathrm{atb},p}^{1,\mathrm{c}}(\M)} \lesssim \|f_\mathrm{d}\|_{\mathrm{h}_{\mathrm{diag}}^1(\M)}.$$ Since the norm in $\mathrm{h}_{\mathrm{diag}}^1(\M)$ is $*$-invariant, we may assume that $f_\mathrm{d}$ is a self-adjoint operator. Then, by an $L_p$-approximation argument we may also assume that the martingale differences have the form $$\Delta_k(f_\mathrm{d}) = \sum_{j \ge 1} \beta_{jk} p_{jk} = \sum_{j \ge 1} \beta_{jk} \Delta_k(p_{jk})$$ for certain $\beta_{jk} \in \R$ and a family $(p_{jk})_{j \ge 1}$ of pairwise disjoint projections. We claim $$|\Delta_k(p)|_{\mathrm{atb},p}^{1,\mathrm{c}} \, \lesssim \, \tau(p)$$ for any projection $p$. This is enough to conclude since then
\begin{eqnarray*}
\|f_\mathrm{d}\|_{\mathrm{H}_{\mathrm{atb},p}^{1,\mathrm{c}}(\M)} & \le & \sum_{j,k \ge 1} |\beta_{jk}| |\Delta_k(p_{jk})|_{\mathrm{atb},p}^{1,\mathrm{c}} \\ & \lesssim & \sum_{j,k \ge 1} |\beta_{jk}|  \tau(p_{jk}) \ = \ \sum_{k \ge 1} \Big\| \sum_{j \ge 1}  \beta_{jk} p_{jk} \Big\|_1 \ = \ \|f_{\mathrm{d}}\|_{\mathrm{h}_{\mathrm{diag}}^1(\M)}. 
\end{eqnarray*}
Let us then prove our claim for $\mathrm{b} = \Delta_k(p)$. To show that $\mathrm{b}$ is a column $p$-atomic block, we start by noticing $\mathsf{E}_{k-1}(\mathrm{b}) = 0$. Let us introduce the family of projections 
\begin{eqnarray*}
q_j(k) & = & \chi_{(\frac{1}{j+1}, \frac{1}{j}]}(\mathsf{E}_k p), \\ q_j(k-1) & = & \chi_{(\frac{1}{j+1}, \frac{1}{j}]}(\mathsf{E}_{k-1} p).
\end{eqnarray*}
Decompose $\mathrm{b}$ into column $p$-subatoms as follows $$\mathrm{b} \, = \, \sum_{j \ge 1} \lambda_j(k) \mathrm{a}_j(k) - \lambda_j(k-1) \mathrm{a}_j(k-1)$$ where coefficients and subatoms are respectively given by 
\begin{eqnarray*}
\lambda_j(k) & = & \mbox{$\frac{2}{j}$} \tau(q_j(k)), \\ \lambda_j(k-1) & = & \mbox{$\frac{1}{j}$} \tau(q_j(k-1)), \\ \mathrm{a}_j(k) & = & \lambda_j(k)^{-1} q_j(k) \mathsf{E}_k(p), \\ \mathrm{a}_j(k-1) & = & \lambda_j(k-1)^{-1} q_j(k-1) \mathsf{E}_{k-1}(p),
\end{eqnarray*}
Since $(q_j(k-1), q_j(k)) \in \M_{k-1} \times \M_k$, we will have a column $p$-atomic block $\mathrm{b}$ if
\begin{itemize}
\item $\mathrm{a}_j(k-1) q_j(k-1) = \mathrm{a}_j(k-1)$ and $\mathrm{a}_j(k) q_j(k) = \mathrm{a}_j(k)$,

\item $\|\mathrm{a}_j(k-1)\|_p \le \tau(q_j(k-1))^{-\frac{1}{p'}}$ and $\|\mathrm{a}_j(k)\|_p \le \frac12 \tau(q_j(k))^{-\frac{1}{p'}}$.
\end{itemize}
It is however a simple exercise to check that this is indeed the case and 
\begin{eqnarray*}
|\mathrm{b}|_{\mathrm{atb},p}^{1,\mathrm{c}} & \le & \sum_{j \ge 1} |\lambda_j(k)| + |\lambda_j(k-1)| \\ & \le & \sum_{j \ge 1}  4 \tau(q_j(k) \mathsf{E}_k(p)) + 2 \tau(q_j(k-1) \mathsf{E}_{k-1}(p)) \ \le \ 6 \tau(p).
\end{eqnarray*}
This justifies our claim above and hence completes the proof of the assertion. \fin

\begin{remark}
\emph{The noncommutative Davis decomposition of Perrin and Junge/Mei \cite{JMe,Pe} is sometimes referred to as the \lq\lq atomic decomposition" for $\mathrm{H}_1^\mathrm{c}(\M)$, since it relates this space with the atomic Hardy space $\mathrm{h}_{\mathrm{at},\mathrm{c}}^1(\M)$ and the diagonal space $\mathrm{h}_{\mathrm{diag}}^1(\M)$. Nevertheless, it seems there is no atomic decomposition of the diagonal part (in the noncommutative setting) beyond the results in this paper.}
\end{remark}

\section{{\bf Two alternative arguments for classical martingales}} \label{Sect1}

In this section we explore two additional proofs of Theorem A valid for classical martingales. None of them work for noncommutative martingales, but shed some light to the problem. The first one uses a weaker notion of atom which proves that atomic blocks can be taken with (at most) two subatoms. Notice that this does not seem to be the case in the von Neumann algebra setting. This is analogous to a similar result for Tolsa's atomic blocks. The second one illustrates how conditional medians instead of medians allow to give a direct proof, avoiding approximation by atomic filtrations. Moreover, we shall obtain in the process an equivalent expression $\|f\|_{\mathrm{BMO}}^\alpha$ for the martingale BMO norm of $f$.   

\subsection{A proof using weak atoms}

Given a probability space $(\Omega, \Sigma, \mu)$ and any filtration $(\Sigma_k)_{k \ge 1}$, we will say that a measurable function $\mathrm{w}: \Omega \to \C$ is a \emph{weak $\infty$-atom} when there is some $k \ge 1$, a $\Sigma_k$-measurable function $\varphi: \Omega \to \C$, with $|\varphi| \le 1$ and $\mathrm{A} := \mathrm{supp} \, \varphi \in \Sigma_k$ so that $$\mathrm{w} = \frac{\varphi - \mathsf{E}_{k-1}(\varphi)}{\mu(\mathrm{A})}.$$ We may find such kind of atoms in \cite{Bo}, but perhaps they were known before. Let us sketch the proof of Theorem A for $p=\infty$ using weak $\infty$-atoms. The proof of the inclusion $\mathrm{H}_{{\mathrm{atb}},\infty}^1(\Omega) \subset \mathrm{H}_1(\Omega)$ will not change from our first proof of Theorem A above. By a straightforward approximation argument (that we will not reproduce here) we may assume that our filtration $(\Sigma_k)_{k \ge 1}$ is atomic. Under this assumption all we need to prove by duality is that $$\|f\|_{\mathrm{BMO}(\Omega)} \, \lesssim \, \sup_{\begin{subarray}{c} \mathrm{b} \, \mathrm{atomic \, block} \\ |\mathrm{b}|_{\mathrm{atb,\infty}}^1 \le 1 \end{subarray}} \, \int_\Omega f \mathrm{b} \, d\mu$$ holds for any $f \in \mathrm{BMO}(\Omega)$. Let us briefly justify this, consider $f \in \mathrm{BMO}(\Omega)$. Since we assume $(\Sigma_k)_{k \ge 1}$ is atomic, given any $\varepsilon > 0$ we may find certain $k \ge 1$ and an atom $\mathrm{A} \in \Sigma_k$ such that $$\| f \|_{\mathrm{BMO}(\Omega)} < (1+\varepsilon) \, \mathsf{E}_k |f - \mathsf{E}_{k-1}(f)|(\mathrm{A}).$$ On the other hand, we have
\begin{eqnarray*}
\lefteqn{\hskip-10pt \mathsf{E}_k |f - \mathsf{E}_{k-1}(f)|(\mathrm{A})} \\ \!\! & = & \!\! \sup_{\|\xi\|_{\infty}\le 1} \frac{1}{\mu(\mathrm{A})} \int_\mathrm{A} \xi (f-\mathsf{E}_{k-1}(f)) \, d\mu \\ \!\! & = & \!\! \sup_{\|\xi\|_{\infty}\le 1} \frac{1}{\mu(\mathrm{A})} \int_\mathrm{A} \xi (f - \mathsf{E}_k(f) + \mathsf{E}_k(f) - \mathsf{E}_{k-1}(f)) \, d\mu \\ \!\! & = & \!\! \sup_{\|\xi\|_{\infty}\le 1} \int_\Omega \underbrace{\chi_{\mathrm{A}}  \frac{\xi - \mathsf{E}_{k}(\xi)}{\mu(\mathrm{A})}}_{\mathrm{a}(\xi)} \, f \, d\mu + \int_\Omega \underbrace{\frac{\mathsf{E}_k(\chi_\mathrm{A} \xi) - \mathsf{E}_{k-1}(\chi_\mathrm{A} \xi)}{\mu(\mathrm{A})}}_{\mathrm{w}(\xi)} \, f \, d\mu. 
\end{eqnarray*}
It is clear that $\mathrm{a}(\xi)$ is an $\infty$-atom with $$|\mathrm{a}(\xi)|_{\mathrm{atb},\infty}^1 \, \le \, 2.$$ Therefore, it suffices to show that $\mathrm{w}(\xi)$ is an $\infty$-atomic block with $$|\mathrm{w}(\xi)|_{\mathrm{atb},\infty}^1 \, \lesssim \, 1.$$ Note that $\mathrm{w}(\xi)$ is a weak $\infty$-atom which can be written as $\mathrm{a}_1(\xi) + \mathrm{a}_2(\xi)$ with
\begin{eqnarray*}
\mathrm{a}_1(\xi) & = & \frac{-1}{\mu(\mathrm{A})} \Big( \frac{1}{\mu(\mathrm{B})} \int_\Omega \chi_\mathrm{A} \xi d\mu \Big) \chi_{\mathrm{B \setminus A}} , \\ \mathrm{a}_2(\xi) & = & \frac{1}{\mu(\mathrm{A})} \Big( \mathsf{E}_k(\xi) - \frac{1}{\mu(\mathrm{B})} \int_\Omega \chi_\mathrm{A} \xi d\mu \Big) \chi_\mathrm{A},
\end{eqnarray*}
where $\mathrm{B}$ is the only atom in $\Sigma_{k-1}$ containing the atom $\mathrm{A} \in \Sigma_k$. Now, it all reduces to show that (up to absolute constants) $\mathrm{a}_1(\xi)$ and $\mathrm{a}_2(\xi)$ are $\infty$-subatoms. Using that $|\xi| \le 1$, we deduce $$|\mathrm{a}_1(\xi)| \, \le \, \frac{\chi_{\mathrm{B \setminus A}}}{\mu(\mathrm{B})} \, \le \, \frac{\chi_{\mathrm{B \setminus A}}}{\mu(\mathrm{B \setminus A})} \quad \mbox{and} \quad |\mathrm{a}_2(\xi)| \, \le \, \frac{2 \chi_\mathrm{A}}{\mu(\mathrm{A})}.$$ Since $\mathrm{a}_1(\xi), \mathrm{a}_2(\xi)$ are $\Sigma_k$-measurable, we easily get the estimate $|\mathrm{w}(\xi)|_{\mathrm{atb},\infty}^1 \le 6$. 

This argument shows that atomic blocks in Theorem A can be taken with at most two subatoms for classical martingales. Unfortunately, the argument does not extend to the noncommutative setting, because of the lack of an approximation argument by atomic filtrations.

\subsection{A proof using conditional medians}

Given a probability space $(\Omega, \Sigma, \mu)$ and a $\sigma$-subalgebra $\Sigma_0\subset \Sigma$, a conditional median $\alpha_0 f$ of a $\Sigma$-measurable function $f$ is a random variable which satisfies:
\begin{itemize}
\item $\alpha_0 f$ is $\Sigma_0$-measurable,
\item Given any $\mathrm{A} \in \Sigma_0$, we have $$\max \Big\{ \mu \big( \mathrm{A} \cap \big\{ f > \alpha_0 f \big\} \big), \mu \big( \mathrm{A} \cap \big\{ f < \alpha_0 f \big\} \big) \Big\} \, \le \, \frac{1}{2} \mu(\mathrm{A}).$$ 
\end{itemize}
Tomkins theorem \cite{Tomkins} shows that each random variable has at least one conditional median with respect to any given $\sigma$-algebra. In the sequel, we will denote a fixed conditional median of $f$ with respect to $\Sigma_k$ by $\alpha_k f$. Before the proof of Theorem A we need a simple lemma which will be crucial in our argument. 

\begin{lemma} \label{CMLemma}
Given $\mathrm{A} \in \Sigma_0$ and $f$ $\Sigma$-measurable $$\mathsf{E}_0 \big( \chi_{\mathrm{A} \cap \{f \le \alpha_0 f\}}\big) \, \ge \, \frac12 \chi_\mathrm{A} \quad \mu\mbox{-a.e.}$$ where $\mathsf{E}_0$ denotes the conditional expectation onto the $\sigma$-subalgebra $\Sigma_0 \subset \Sigma$. 
\end{lemma}

\dem By the definition of conditional median $$\mu \big( \mathrm{B} \cap \big\{ f \le \alpha_0f \big\} \big) \ge \frac12 \mu(\mathrm{B})$$ for every $\Sigma_0$-measurable set $\mathrm{B}$. Assume now that the set $\mathrm{A}$ in the statement fails the given inequality and define $\mathrm{B}$ to be the $\Sigma_0$-measurable level set where $\mathsf{E}_0 (\chi_{\mathrm{A} \cap \{f \le \alpha_0 f\}}) < \frac12$. If the assertion failed for $\mathrm{A}$, we would have $\mu(\mathrm{B}) > 0$ and we could conclude that 
\begin{eqnarray*}
\mu \big( \mathrm{B} \cap \big\{ f \le \alpha_0f \big\} \big) & = & \int_\mathrm{B} \mathsf{E}_0 \big( \chi_{\mathrm{B} \cap \{f \le \alpha_0 f\}}\big) \, d\mu \\ & \le & \int_\mathrm{B} \mathsf{E}_0 \big( \chi_{\mathrm{A} \cap \{f \le \alpha_0 f\}}\big) \, d\mu \ < \ \frac12 \mu(\mathrm{B})
\end{eqnarray*}
which contradicts the definition of conditional median. The proof is complete. \fin

\demA
Again, the proof of $\mathrm{H}_{{\mathrm{atb}},p}^1(\Omega) \subset \mathrm{H}_1(\Omega)$ will not change from our first proof of Theorem A above. Here we only prove the reverse inclusion $\mathrm{H}_1(\Omega) \subset \mathrm{H}_{{\mathrm{atb}},p}^1(\Omega)$. To that end, we will show that $$\mathrm{H}_{{\mathrm{atb}},p}^1(\Omega)^* \subset \mathrm{BMO}(\Omega),$$ which suffices by duality. Let $L: \mathrm{H}_{{\mathrm{atb}},p}^1(\Omega) \to \C$ be a continuous functional in the dual space. To proceed, we need to show that $L = L_f$ acts by integration in $(\Omega,\mu)$ against a function $f \in L_{\mathrm{loc}}^1(\Omega)$ and deduce a posteriori that $f \in \mathrm{BMO}(\Omega)$ and we have $$\|f\|_{\mathrm{BMO}(\Omega)} \, \le \, \mathrm{C}_p \|L_f\|_{\mathrm{H}_{\mathrm{atb},p}^1(\Omega)^*}$$ for some absolute constant $\mathrm{C}_p$. The existence of such $f$ follows from the inclusion $\mathrm{h}_{\mathrm{at},p}^1(\Omega) \subset \mathrm{H}_{\mathrm{atb},p}^1(\Omega)$, so that $\mathrm{H}_{\mathrm{atb},p}^1(\Omega)^* \subset \mathrm{h}_{\mathrm{at},p}^1(\Omega)^* = \mathrm{bmo}(\Omega)$. In particular any continuous functional $L$ in the dual of $\mathrm{H}_{\mathrm{atb},p}^1(\Omega)$ can be represented by a function $f \in \mathrm{bmo}(\Omega)$. We now claim that $$\frac{1}{\mathrm{c}_p} \|f\|_{\mathrm{BMO}(\Omega)} \, \le \, \|f\|_{\mathrm{BMO}}^\alpha \, \le \, \mathrm{c}_p \|L_f\|_{\mathrm{H}_{\mathrm{atb},p}^1(\Omega)^*},$$ where $\mathrm{c}_p$ only depends on $p$ and $\|f\|_{\mathrm{BMO}}^\alpha$ is given by $$\|f\|_{\mathrm{BMO}}^\alpha \, = \, \max \Big\{ \|\mathsf{E}_1f\|_\infty, \ \sup_{k \ge 1} \big\| \mathsf{E}_k | f - \alpha_k f |^{p'} \big\|_\infty^{\frac{1}{p'}}, \ \sup_{k \ge 2} \big\| \alpha_{k}f - \alpha_{k-1} f \big\|_\infty \Big\}.$$ Note that this quantity depends a priori on the choice of the conditional medians $\alpha_k f$. This however will be unsubstantial since our inequalities hold with constants which are independent of our choice. It is clear that the proof will be complete if we justify our claim, which we will in two steps.

\noindent \textbf{Step 1.} The inequality $$\|f\|_{\mathrm{BMO}(\Omega)} \le \mathrm{c}_p \|f\|_{\mathrm{BMO}^\alpha}$$ is the simplest one. Namely, by John-Nirenberg inequality we have
\begin{eqnarray*}
\|f\|_{\mathrm{BMO}(\Omega)} & = & \sup_{k \ge 1} \big\| \mathsf{E}_k |f - \mathsf{E}_{k-1}f|^2 \big\|_\infty^{\frac12} \\ & \sim & \|\mathsf{E}_1 f\|_\infty + \sup_{k \ge 1} \big\| \mathsf{E}_k |f - \mathsf{E}_k f |^2 \big\|_\infty^{\frac12} + \sup_{k \ge 2} \|df_k\|_\infty \\ & \sim & \|\mathsf{E}_1 f\|_\infty + \sup_{k \ge 1} \big\| \mathsf{E}_k |f - \mathsf{E}_k f |^{p'} \big\|_\infty^{\frac1{p'}} + \sup_{k \ge 2} \|df_k\|_\infty \ = \ \mathrm{A}_1 + \mathrm{A}_2 + \mathrm{A}_3.
\end{eqnarray*}
The term $\mathrm{A}_1$ admits a trivial bound. Next 
\begin{eqnarray*}
\mathrm{A}_2 & \le & \sup_{k \ge 1} \big\| \mathsf{E}_k |f - \alpha_k f |^{p'} \big\|_\infty^{\frac1{p'}} + \big\| \mathsf{E}_k |\alpha_k f - \mathsf{E}_k f |^{p'} \big\|_\infty^{\frac1{p'}} \\ & \le & \|f\|_{\mathrm{BMO}}^\alpha + \sup_{k \ge 1} \big\| \mathsf{E}_k (f - \alpha_k f) \big\|_\infty \ \le \ 2 \, \|f\|_{\mathrm{BMO}}^\alpha,
\end{eqnarray*}
where the last inequality uses conditional Jensen's inequality $\phi(\mathsf{E}_k f) \le \mathsf{E}_k(\phi(f))$ for the convex function $\phi(x) = x^{p'}$. Finally, the last term $\mathrm{A}_3$ is estimated by decomposing $df_k = \mathsf{E}_k(f - \alpha_k f) + (\alpha_k f - \alpha_{k-1} f) - \mathsf{E}_{k-1}(f - \alpha_{k-1}f) $ together with the triangle inequality and conditional Jensen's inequality one more time.  

\noindent \textbf{Step 2.} The inequality $$\|f\|_{\mathrm{BMO}}^{\alpha} \le \mathrm{c}_p \|L_f\|_{\mathrm{H}_{\mathrm{atb},p}^1(\Omega)^*}$$ requires a bit more work. Since $\Sigma_1$-measurable functions are atomic blocks $$ \|\mathsf{E}_1f\|_{\infty} = \sup_{\mathrm{B} \in \Sigma_1} \Big| \mean_\mathrm{B} f \, d\mu \Big| \le \frac{1}{\mu(\mathrm{B})} \|L_f\|_{\mathrm{H}_{\mathrm{atb},p}^1(\Omega)^*} \|\chi_\mathrm{B}\|_{\mathrm{H}_{\mathrm{atb},p}^1(\Omega)} \le \|L_f\|_{\mathrm{H}_{\mathrm{atb},p}^1(\Omega)^*}.$$ Let us now bound the other two terms in $\|f\|_{\mathrm{BMO}}^\alpha$. In order to estimate the second term, we will use that for any $\mathrm{A} \in \Sigma_k$ there exists a $p$-atomic block $\mathrm{b}_{\mathrm{A},f} \neq 0$ satisfying the following two inequalities $$\|\mathrm{b}_{\mathrm{A},f}\|_{\mathrm{H}_{\mathrm{atb},p}^1(\Omega)} \, \lesssim \, \mu(\mathrm{A})^{\frac{1}{p'}} \Big( \int_\mathrm{A} |f - \alpha_k f|^{p'} \, d\mu \Big)^{\frac{1}{p}} \, \lesssim \, \mu(\mathrm{A})^{\frac{1}{p'}} \Big| \int_\Omega f \mathrm{b}_{\mathrm{A},f} \, d\mu \Big|^{\frac{1}{p}}.$$ This immediately implies that $$\sup_{k \ge 1} \big\| \mathsf{E}_k | f - \alpha_k f |^{p'} \big\|_\infty^{\frac{1}{p'}} \, \lesssim \, \|L_f\|_{\mathrm{H}_{\mathrm{atb},p}^1(\Omega)^*}$$ as desired. Indeed, this can be justified as follows
\begin{eqnarray*}
\big\| \mathsf{E}_k | f - \alpha_k f |^{p'} \big\|_\infty^{\frac{1}{p'}} & = & \sup_{\mathrm{A} \in \Sigma_k} \Big( \mean_\mathrm{A} | f - \alpha_k f |^{p'} \, d\mu \Big)^{\frac{1}{p'}} \\ & \lesssim & \sup_{\mathrm{A} \in \Sigma_k} \frac{1}{\|\mathrm{b}_{\mathrm{A},f}\|_{\mathrm{H}_{\mathrm{atb},p}^1(\Omega)}} \Big| \int_\Omega f \mathrm{b}_{\mathrm{A},f} \, d\mu \Big| \ \le \ \|L_f\|_{\mathrm{H}_{\mathrm{atb},p}^1(\Omega)^*}.
\end{eqnarray*}
Given $\mathrm{A} \in \Sigma_k$, let us then prove the existence of such $p$-atomic block. Assume $$\int_{\mathrm{A} \cap \{f>\alpha_kf\}} |f-\alpha_k f|^{p'} \, d\mu \, \ge \, \int_{\mathrm{A} \cap \{f<\alpha_kf\}} |f-\alpha_kf|^{p'} \, d\mu.$$ This assumption is admissible since we may easily modify the construction of our $p$-atomic block $\mathrm{b}_{\mathrm{A},f}$ to satisfy the required estimates in case the inequality above is reversed. Define the function $$\mathrm{b}_{\mathrm{A},f}(x) \, = \,  
|f-\alpha_kf|^{p'-1} \chi_{\mathrm{A} \cap \{f>\alpha_k f\}} - \frac{\mathsf{E}_k(|f-\alpha_kf|^{p'-1} \chi_{\mathrm{A} \cap \{f>\alpha_k f\}})}{\mathsf{E}_k(\chi_{\mathrm{A} \cap \{f \le \alpha_kf\}})} \chi_{\mathrm{A} \cap \{f \le \alpha_k f\}}.$$
Obviously, $\mathsf{E}_k(\mathrm{b}_{\mathrm{A},f})=0$ and $\mathrm{supp} (\mathrm{b}_{\mathrm{A},f}) \subset \mathrm{A}$. This yields
\begin{eqnarray*} 
\|\mathrm{b}_{\mathrm{A},f}\|_{\mathrm{H}_{\mathrm{atb},p}^1(\Omega)} & \le & \mu(\mathrm{A})^{\frac{1}{p'}} \|\mathrm{b}_{\mathrm{A},f}\|_p \\ [6pt] & \le & \mu(\mathrm{A})^{\frac{1}{p'}} \Big( \int_{\mathrm{A} \cap\{f>\alpha_kf\}} |f-\alpha_kf|^{p(p'-1)} d\mu \Big)^{\frac{1}{p}} \\ & + & \mu(\mathrm{A})^{\frac{1}{p'}} \Big( \int_{\mathrm{A} \cap\{f\leq \alpha_kf\}} \Big[ \frac{\mathsf{E}_k(|f-\alpha_kf|^{p'-1})}{[\mathsf{E}_k(\chi_{\mathrm{A} \cap\{f \leq \alpha_kf\}})]} \Big]^p d\mu \Big)^{\frac{1}{p}} \ = \ \mathrm{A}_1 + \mathrm{A}_2.
\end{eqnarray*}
Since $p(p'-1)=p'$, $\mathrm{A}_1$ clearly satisfies the desired estimate. On the other hand 
\begin{eqnarray*}
\mathrm{A}_2 \!\!\! & = & \!\!\! \mu(\mathrm{A})^{\frac{1}{p'}} \Big( \int_\mathrm{A} \chi_{\mathrm{A} \cap \{f \leq \alpha_kf\}} \Big[ \frac{\mathsf{E}_k(|f-\alpha_kf|^{p'-1})}{[\mathsf{E}_k(\chi_{\mathrm{A} \cap\{f \leq \alpha_kf\}})]} \Big]^p \, d\mu \Big)^{\frac1p} \\ 
\!\!\! & = & \!\!\! \mu(\mathrm{A})^{\frac{1}{p'}} \Big( \int_\mathrm{A} [\mathsf{E}_k(|f-\alpha_kf|^{p'-1})]^p [\mathsf{E}_k(\chi_{\mathrm{A} \cap\{f \leq \alpha_kf\}})]^{1-p} \, d\mu \Big)^{\frac1p} \\ \!\!\! & \lesssim & \!\!\! \mu(\mathrm{A})^{\frac{1}{p'}} \Big(  \int_\mathrm{A} [\mathsf{E}_k(|f-\alpha_kf|^{p'-1})]^p \, d\mu \Big)^{\frac1p} \, \le \, \mu(\mathrm{A})^{\frac{1}{p'}} \Big(  \int_\mathrm{A} \mathsf{E}_k(|f-\alpha_kf|^{p'}) \, d\mu \Big)^{\frac1p},
\end{eqnarray*}
where we have used Lemma \ref{CMLemma} for the first inequality and conditional Jensen's inequality for the second one. Now, since $\mathrm{A} \in \Sigma_k$, we can remove the conditional expectation $\mathsf{E}_k$ in the integrand of the last term above 
to complete the proof of the estimate for $|\mathrm{b}_{\mathrm{A},f}|_{\mathrm{atb},p}^1$. The other inequality is simpler. Since $(f-\alpha_kf) \mathrm{b}_{\mathrm{A},f}$ is nonnegative by definition of $\mathrm{b}_{\mathrm{A},f}$ and $\mathsf{E}_k(\mathrm{b}_{\mathrm{A},f}) = 0$, we get
\begin{eqnarray*}
\int_\Omega f \mathrm{b}_{\mathrm{A},f} \, d\mu & = & \int_\Omega (f-\alpha_k f) \mathrm{b}_{\mathrm{A},f} \, d\mu \\
& \ge & \int_{\mathrm{A} \cap \{f > \alpha_k f\}} |f-\alpha_k f|^{p'} \, d\mu \ \ge \ \frac{1}{2} \int_{\mathrm{A}} |f-\alpha_k f|^{p'} \, d\mu.
\end{eqnarray*}
This completes the proof of the expected estimate for the second term in $\|f\|_{\mathrm{BMO}}^\alpha$. It remains to prove that $$\sup_{k \ge 2} \big\| \alpha_{k}f - \alpha_{k-1} f \big\|_\infty \, = \, \sup_{k \ge 2} \sup_{\mathrm{A} \in \Sigma_k} \mean_{\mathrm{A}} |\alpha_k f - \alpha_{k-1} f | \, d\mu  \, \lesssim \, \mathrm{c}_p \|L_f\|_{\mathrm{H}_{\mathrm{atb},p}^1(\Omega)^*}.$$ Fix $k>1$ and $\mathrm{A} \in \Sigma_k$. By the triangle and Jensen's inequality $$\mean_{\mathrm{A}} |\alpha_kf - \alpha_{k-1}f| \, d\mu \, \le \, \Big( \mean_\mathrm{A} |f - \alpha_{k}f|^{p'} \, d\mu \Big)^{\frac{1}{p'}} + \Big( \mean_\mathrm{A} | f - \alpha_{k-1}f |^{p'} d\mu \Big)^{\frac{1}{p'}}.$$ Since $\mathrm{A} \in \Sigma_k$, the first term in the right hand side is bounded above by $$\big\| \mathsf{E}_k \big| f - \alpha_k f \big|^{p'} \big\|_\infty^{\frac{1}{p'}} \, \lesssim \, \|L_f\|_{\mathrm{H}_{\mathrm{atb},p}^1(\Omega)^*}$$ as we proved before. To bound the second term, we consider the function $$\mathrm{b}_{\mathrm{A},f} = \underbrace{\frac{|f-\alpha_{k-1}f|^{p'}}{f-\alpha_{k-1}f}\chi_{\mathrm{A} \cap\{f\not= \alpha_{k-1}f\}}}_{\lambda_* \mathrm{a}_*} - \underbrace{\mathsf{E}_{k-1} \Big( \frac{|f-\alpha_{k-1}f|^{p'}}{f-\alpha_{k-1}f} \chi_{\mathrm{A} \cap\{f\not= \alpha_{k-1}f\}} \Big)}_{\sum_{j \in \Z} \lambda_j \mathrm{a}_j}$$ where $$\lambda_j \mathrm{a}_j = \mathsf{E}_{k-1} \Big( \frac{|f-\alpha_{k-1}f|^{p'}}{f-\alpha_{k-1}f} \chi_{\mathrm{A} \cap\{f\not= \alpha_{k-1}f\}} \Big) \underbrace{\chi_{\{2^{j-1} < \mathsf{E}_{k-1}(|f - \alpha_{k-1}f|^{p'-1} \chi_{\mathrm{A}}) \le 2^j\}}}_{\chi_{\mathrm{B}_j}}.$$
We have $\mathsf{E}_{k-1}(\mathrm{b}_{\mathrm{A},f}) = 0$ so that 
\begin{eqnarray*}
\|\mathrm{b}_{\mathrm{A},f}\|_{\mathrm{H}_{\mathrm{atb},p}^1(\Omega)} & \le & |\lambda_*| + \sum_{j \in \Z} |\lambda_j| \\ & \le & \mu(\mathrm{A})^{\frac{1}{p'}} \big\| |f-\alpha_{k-1}f|^{p'-1} \chi_\mathrm{A} \big\|_p \\ [6pt] & + & \sum_{j \in \Z} \mu(\mathrm{B}_j)^{\frac{1}{p'}} \big\| \mathsf{E}_{k-1} \big( |f - \alpha_{k-1}f|^{p'-1} \chi_{\mathrm{A}} \big) \chi_{\mathrm{B}_j} \big\|_p
\end{eqnarray*}
The second term in the right hand side is dominated by the first one since
\begin{eqnarray*}
\lefteqn{\hskip-30pt \sum_{j \in \Z} \mu(\mathrm{B}_j)^{\frac{1}{p'}} \big\| \mathsf{E}_{k-1} \big( |f - \alpha_{k-1}f|^{p'-1} \chi_{\mathrm{A}} \big) \chi_{\mathrm{B}_j} \big\|_p} \\ & \le &  \sum_{j \in \Z} 2^j \mu(\mathrm{B}_j) \\ & \sim & \sum_{j \in \Z} \mean_{\mathrm{B}_j} \mathsf{E}_{k-1} \big( |f-\alpha_{k-1}f|^{p'-1} \chi_{\mathrm{A}} \big) \, d\mu \ \mu(\mathrm{B}_j) \\ & = & \int_{\cup \mathrm{B}_j} \mathsf{E}_{k-1} \big( |f-\alpha_{k-1}f|^{p'-1} \chi_{\mathrm{A}} \big) d\mu \\ & = & \int_\mathrm{A} |f-\alpha_{k-1}f|^{p'-1} d\mu \ \le \ \mu(\mathrm{A})^{\frac{1}{p'}} \big\| |f-\alpha_{k-1}f|^{p'-1}\chi_\mathrm{A} \big\|_p.
\end{eqnarray*}
In summary, we have proved that $$\|\mathrm{b}_{\mathrm{A},f}\|_{\mathrm{H}_{\mathrm{atb},p}^1(\Omega)} \, \lesssim \, \mu(\mathrm{A})^{\frac{1}{p'}} \big\| |f-\alpha_{k-1}f|^{p'-1}\chi_\mathrm{A} \big\|_p.$$ On the other hand, let us observe that
\begin{eqnarray*}
\int_\Omega f \mathrm{b}_{\mathrm{A},f} \, d\mu & = & \int_\Omega (f-\alpha_{k-1}f) \mathrm{b}_{\mathrm{A},f} \, d\mu \\ & = & \int_\mathrm{A} |f-\alpha_{k-1}f|^{p'} \, d\mu \\ & - & \int_\Omega (f-\alpha_{k-1}f) \mathsf{E}_{k-1} \Big( \frac{|f-\alpha_{k-1}f|^{p'}}{f-\alpha_{k-1}f} \chi_{A\cap\{f\not= \alpha_{k-1}f\}} \Big) \, d\mu.\\
\end{eqnarray*}
Using this and the estimates so far we obtain 
\begin{eqnarray*}
\lefteqn{\hskip-10pt \int_\mathrm{A} |f-\alpha_{k-1}f|^{p'} \, d\mu \ \le \ \Big| \int_\Omega f \mathrm{b}_{\mathrm{A},f} \, d\mu \Big|} \\ & + & \summ_j \int_{\mathrm{B}_j} |f-\alpha_{k-1}f| \, \mathsf{E}_{k-1} \big( |f-\alpha_{k-1}f|^{p'-1} \chi_{\mathrm{A}} \big) \, d\mu \\ [1pt] & \le & \|L_f\|_{\mathrm{H}_{\mathrm{atb},p}^1(\Omega)^*} \|\mathrm{b}_{\mathrm{A},f}\|_{\mathrm{H}_{\mathrm{atb},p}^1(\Omega)} \\ [3pt] & + & \summ_j \big\| \mathsf{E}_{k-1} |f-\alpha_{k-1}f|^{p'} \big\|_\infty^{\frac{1}{p'}} \mu(\mathrm{B}_j)^{\frac{1}{p'}} \big\| \mathsf{E}_{k-1} \big( |f-\alpha_{k-1}f|^{p'-1} \chi_{\mathrm{A}} \big) \chi_{\mathrm{B}_j} \big\|_{p} \\ [4pt] & \lesssim & \|L_f\|_{\mathrm{H}_{\mathrm{atb},p}^1(\Omega)^*} \|\mathrm{b}_{\mathrm{A},f}\|_{\mathrm{H}_{\mathrm{atb},p}^1(\Omega)} \\
[1pt] & + & \|L_f\|_{\mathrm{H}_{\mathrm{atb},p}^1(\Omega)^*} \summ_j \mu(\mathrm{B}_j)^{\frac{1}{p'}} \big\| \mathsf{E}_{k-1} \big( |f-\alpha_{k-1}f|^{p'-1} \chi_{\mathrm{A}} \big) \chi_{\mathrm{B}_j} \big\|_{p} \\ & \lesssim & \|L_f\|_{\mathrm{H}_{\mathrm{atb},p}^1(\Omega)^*} \mu(\mathrm{A})^{\frac{1}{p'}} \big\| |f-\alpha_{k-1}f|^{p'-1} \chi_\mathrm{A} \big\|_p.
\end{eqnarray*}
Rearranging and noticing that $p(p'-1)=p'$ we get $$\Big( \mean_\mathrm{A} |f-\alpha_{k-1}f|^{p'} \, d\mu \Big)^{\frac{1}{p'}} \, \lesssim \, \|L_f\|_{\mathrm{H}_{\mathrm{atb},p}^1(\Omega)^*},$$
the desired estimate. This completes the proof of Theorem A for $p < \infty$. \fin

\demAA
The proof presents a lot of similarities with the case $p < \infty$. As above, we will only prove the inclusion $\mathrm{H}_1(\Omega) \subset \mathrm{H}_{\mathrm{atb},\infty}^1(\Omega)$. Again, we proceed by duality and the goal is to show that $$\|f\|_{\mathrm{BMO}(\Omega)} \, \lesssim \, \|f\|_{\mathrm{BMO}}^\alpha \, \lesssim \, \|L_f\|_{\mathrm{H}_{\mathrm{atb},\infty}^1(\Omega)^*}.$$ Our former argument for the first inequality is still valid. Now consider 
\begin{enumerate}
\item There exists $k \ge 1$ and $\mathrm{A} \in \Sigma_k$ such that $$\mean_\mathrm{A} |f-\alpha_kf| \, d\mu \, \ge \, \frac{1}{32} \|f\|_{\mathrm{BMO}}^{\alpha}.$$

\item Property (1) fails and there exists $k \ge 2$ such that $$\big\| \alpha_{k} f - \alpha_{k-1} f \big\|_\infty \, \ge \, \frac{1}{2} \|f\|_{\mathrm{BMO}}^{\alpha}.$$

\item The following inequality holds $$\hskip10pt \max \Big\{ \sup_{k \ge 1} \big\| \mathsf{E}_k | f - \alpha_k f |^{p'} \big\|_\infty^{\frac{1}{p'}}, \ \sup_{k \ge 2} \big\| \alpha_{k}f - \alpha_{k-1} f \big\|_\infty \Big\} \, \le \, \|\mathsf{E}_1f\|_\infty.$$ 
\end{enumerate}
It is not difficult to check that at least one of the properties above always hold for every $f$ with $\|f\|_{\mathrm{BMO}}^\alpha$ finite. When (3) holds, we may argue as in the proof of the case $p < \infty$ to deduce $\|f\|_{\mathrm{BMO}}^{\alpha} \le \|L_f\|_{(\mathrm{H}_{\mathrm{atb},\infty}^1(\Omega))^*}$. When $(1)$ holds, we consider the following function $$\mathrm{b}_{\mathrm{A},f} = \underbrace{\chi_{\mathrm{A} \cap \{f>\alpha_k f\}\}} - \chi_{\mathrm{A} \cap \{f<\alpha_k f\}}}_{\mathrm{b}_1} -\underbrace{\chi_{\mathrm{A} \cap \{f=\alpha_kf\}} \mathsf{E}_k(\mathrm{b}_1) [\mathsf{E}_k(\chi_{\mathrm{A} \cap \{f=\alpha_kf\}}]^{-1}}_{\mathrm{b}_2},$$ with the convention $0.\infty = 0$ when $\mathrm{A} \cap \{f = \alpha_k f\} = \emptyset$. Obviously, $
\mathsf{E}_k(\mathrm{b}_{\mathrm{A},f})=0$ and $\|\mathrm{b}_1\|_{\infty} \le 1$. Decomposing into level sets as we did in the proof for $p < \infty$, one can show that $\|\mathrm{b}_2\|_\infty \le 4$, details are 
left to the reader. These $L_\infty$ estimates yield $$\|\mathrm{b}_{\mathrm{A},f}\|_{\mathrm{H}_{\mathrm{atb},\infty}^1(\Omega)} \, \lesssim \, \mu(\mathrm{A}).$$ Moreover, we have $$\Big| \int_\Omega f \mathrm{b}_{\mathrm{A},f} \, d\mu \Big| = \Big| \int_\Omega (f-\alpha_k f) \mathrm{b}_1 \, d\mu \Big| = \int_\mathrm{A} |f-\alpha_kf| \, d\mu \ge \frac{1}{32} \|f\|_{\mathrm{BMO}}^{\alpha} \mu(\mathrm{A})$$ by assumption $(1)$. This implies $$\|L_f\|_{\mathrm{H}_{\mathrm{atb},\infty}^1(\Omega)^*} \|\mathrm{b}_{\mathrm{A},f}\|_{\mathrm{H}_{\mathrm{atb},\infty}^1(\Omega)} \ge \frac{1}{32} \|f\|_{\rm BMO}^{\alpha} \mu(\mathrm{A}) \gtrsim \frac{1}{32} \|f\|_{\mathrm{BMO}}^{\alpha} \|\mathrm{b}_{\mathrm{A},f}\|_{\mathrm{H}_{\mathrm{atb},\infty}^1(\Omega)},$$ which is what we wanted. Finally, if $(2)$ holds there exists $\mathrm{A} \in \Sigma_k$ such that $$\Big| \mean_\mathrm{A} ( \alpha_k f -\alpha_{k-1}f ) \, d\mu \Big| > \frac{1}{4} \|f\|_{\mathrm{BMO}}^{\alpha}.$$ Let $\mathrm{B} = \mathrm{supp}(\mathsf{E}_{k-1}(\chi_\mathrm{A})) \in \Sigma_{k-1}$. Define $\mathrm{b}_{\mathrm{A},f}$ in this case as $$\mathrm{b}_{\mathrm{A},f} \, = \, \chi_\mathrm{A} - \mathsf{E}_{k-1}(\chi_\mathrm{A}).$$ Obviously, it is a $\infty$-atomic block. Taking $\mathrm{B}_j = \{(j-1)/\mathrm{N} < \mathsf{E}_{k-1}(\chi_\mathrm{A}) \leq j/\mathrm{N}\}$, we see that $$\|\mathrm{b}_{\mathrm{A},f}\|_{\mathrm{H}_{\mathrm{atb}, \infty}^1(\Omega)} \, \lesssim \, \mu(\mathrm{A}) + \sum_{j=1}^\mathrm{N} \big\| \mathsf{E}_{k-1}(\chi_\mathrm{A}) \chi_{\mathrm{B}_j} \big\|_{\infty} \mu(\mathrm{B}_j),$$ for all $\mathrm{N}$. The sum in the right hand side converges to $$\int_\Omega \mathsf{E}_{k-1}(\chi_\mathrm{A}) \, d\mu \, = \, \mu(\mathrm{A})$$ as $\mathrm{N} \rightarrow \infty$. This shows that $\|\mathrm{b}_{\mathrm{A},f}\|_{\mathrm{H}_{\mathrm{atb}, \infty}^1(\Omega)} \, \lesssim \, \mu(\mathrm{A})$. Next we compute 
\begin{eqnarray*}
L_f(\mathrm{b}_{\mathrm{A},f}) & = & \int_{\mathrm{B}} \mathrm{b}_{\mathrm{A},f} (f-\alpha_{k-1}f) \, d\mu \\ & = & \int_\mathrm{A} (f-\alpha_{k-1}f) d\mu - \int_{\mathrm{B}} \mathsf{E}_{k-1}(\chi_A)(f-\alpha_{k-1}f) \, d\mu \\ & = & \int_\mathrm{A} (f-\alpha_{k}f) \, d\mu + \int_\mathrm{A} (\alpha_{k}f - \alpha_{k-1}f) \, d\mu - \int_\mathrm{B} \mathsf{E}_{k-1}(\chi_\mathrm{A})(f-\alpha_{k-1}f) \, d\mu.
\end{eqnarray*}
Since $(1)$ does not hold, we have $$\Big| \int_\mathrm{A} (f-\alpha_{k}f) \, d\mu \Big| \, \le \, \frac{1}{32} \|f\|_{\mathrm{BMO}}^{\alpha} \mu(\mathrm{A}).$$
On the other hand, and splitting into level sets we find 
\begin{eqnarray*}
\Big| \int_\mathrm{B} \mathsf{E}_{k-1}(\chi_\mathrm{A})(f-\alpha_{k-1}f) \, d\mu \Big| & \le & \sum_{j=1}^\mathrm{N} \frac{j}{\mathrm{N}} \Big| \int_{\mathrm{B}_j} f-\alpha_{k-1}f \, d\mu \Big| \\ & = & \sum_{j=1}^\mathrm{N} \frac{j}{\mathrm{N}} \mu(\mathrm{B}_j) \Big| \mean_{\mathrm{B}_j} f-\alpha_{k-1}f \, d\mu \Big| \\ & \leq & \sup_{\mathrm{C} \in \Sigma_{k-1}} \Big| \mean_\mathrm{C} (f-\alpha_{k-1}f) \, d\mu \Big| \sum_{j=1}^\mathrm{N} \frac{j}{\mathrm{N}} \mu(\mathrm{B}_j)
\end{eqnarray*}
which is dominated by $\frac{1}{16} \|f\|_{\mathrm{BMO}}^{\alpha} \mu(\mathrm{A})$ for $\mathrm{N}$ large enough. So we get
\begin{eqnarray*}
\|L_f\|_{\mathrm{H}_{\mathrm{atb},\infty}^1(\Omega)^*} & \ge & \frac{1}{\|\mathrm{b}_{\mathrm{A},f}\|_{\mathrm{H}_{\mathrm{atb},\infty}^1(\Omega)}} |L_f(\mathrm{b}_{\mathrm{A},f})| \\ & \gtrsim & \frac{1}{\mu(\mathrm{A})} \big( \frac12 - \frac{1}{32} - \frac{1}{16} \big) \|f\|_{\mathrm{BMO}}^\alpha \mu(\mathrm{A}) \ \gtrsim \ \|f\|_{\mathrm{BMO}}^\alpha.
\end{eqnarray*}
This is the last possible case and completes the proof of Theorem A for $p = \infty$. \fin

\section{{\bf Open problems}}

When $0<p<1$, one can extend the definition of atomic blocks to $(p,q)$-atomic blocks. Given $1<q< \infty$, $\mathrm{b}$ is called a $(p,q)$-atomic block when $\mathrm{b}$ is $\Sigma_1$-measurable or there exists $k \ge 1$ such that the following properties hold:
\begin{itemize}
\item $\mathsf{E}_k(\mathrm{b}) = 0$,

\vskip2pt

\item $\mathrm{b} = \sum_j \lambda_j \mathrm{a}_j$ where 

\vskip4pt

\begin{itemize}
\item[$\circ$] $\mathrm{supp}(\mathrm{a}_j) \subset \mathrm{A}_j$, 

\item[$\circ$] $\|\mathrm{a}_j\|_q \le \mu(\mathrm{A}_j)^{1-\frac{1}{p}-\frac{1}{q'}} \frac{1}{k_j - k + 1}$,
\end{itemize}

\vskip3pt

\noindent for certain $k_j \ge k$ and $\mathrm{A}_j \in \Sigma_{k_j}$. 
\end{itemize}
As in the case of $p=1$, set $|\mathrm{b}|_{\mathrm{atb},q}^p = \|\mathrm{b}\|_p$ if $\mathrm{b} \in L_p(\Omega, \Sigma_1,\mu)$ and $$|\mathrm{b}|_{\mathrm{atb},q}^p \, = \, \inf_{\begin{subarray}{c} \mathrm{b} = \sum_j \lambda_j \mathrm{a}_j \\ \mathrm{a}_j \ (p,q)-\mathrm{subatom} \end{subarray}} \ \sum_{j \ge 1} |\lambda_j|$$ otherwise. Finally, we define 
$$\mathrm{H}_{\mathrm{atb},q}^p(\Omega) \, = \, \Big\{ f \in L_p(\Omega) \ \big| \ f = \summ_i \mathrm{b}_i, \ \mathrm{b}_i \ \mbox{martingale } \mbox{$(p,q)$-atomic block} \Big\},$$
equipped with the following quasi-norm $$\|f\|_{\mathrm{H}_{\mathrm{atb},q}^p(\Omega)} \, = \, \inf_{\begin{subarray}{c} f = \sum_i \mathrm{b}_i \\ \mathrm{b}_i \ (p,q)-\mathrm{atomic \, block} \end{subarray}} \Big[ \sum_{i \ge 1} \big( |\mathrm{b}_i|_{\mathrm{atb},q}^p \big)^p \Big]^{\frac{1}{p}}.$$ The spaces $\mathrm{H}_{\mathrm{atb},q}^p(\Omega)$ defined above are quasi-Banach subspaces of $L_p(\Omega)$. One can follow almost \emph{verbatim} the steps in the proof of Theorem A to conclude that the set of linear continuous functionals acting on $\mathrm{H}_{\mathrm{atb},q}^p(\Omega)$ can be identified with the Lipschitz type class $\Lambda_{p,q}(\Omega)$ of functions with finite norm $$\|f\|_{\Lambda_{p,q}(\Omega)} \, = \, \sup_{\begin{subarray}{c} k \ge 1 \\ \mathrm{A} \in \Sigma_k \end{subarray}} \frac{1}{\mu(\mathrm{A})^{\frac{1}{p} -1}} \Big[ \left(\mean_A |f-\mathsf{E}_k f|^q \ d\mu \right)^{\frac{1}{q}} + \left\| df_k \right\|_{\infty} \Big].$$ Notice that when $p \rightarrow 1$, the norm in $\Lambda_{p,q}(\Omega)$ tends to the norm in ${\rm BMO}(\Omega)$. This motivates our first problem, which is somehow analogous (in the context of atomic blocks of this paper) to Problem $3$ in \cite{BCP}.

\begin{problem}
Do we have $$\mathrm{H}_{{\rm atb},q}^p(\Omega) = \mathrm{H}_{p} (\Omega)$$ for $0<p<1$ and $q>1$\emph{?} Moreover, do we have $\mathrm{H}_{p} (\Omega)^* = \Lambda_{p,q}(\Omega)$\emph{?}
\end{problem}

Our main result shows that a function in $\mathrm{H}_{1}(\Omega)$ can be decomposed into atomic blocks similar to the ones appearing in the definition of $\mathrm{H}_{\rm atb}^1 (\mathbb{R}^n,\mu)$, the atomic block Hardy space of Tolsa \cite{To}. In the proof given in Section 1, we make use of Davis decomposition for martingales. It is natural to ask whether we can find a description of the space $\mathrm{H}_{\rm atb}^1 (\mathbb{R}^n,\mu)$ in terms of some sort of Davis decomposition that splits the space into a (classical) atomic part and a diagonal part. Note that a suitable candidate for the atomic part is the space $\mathrm{h}_{\rm at}^1 (\mathbb{R}^n,\mu)$ of functions decomposable into classical atoms supported on doubling sets, since in that case one can easily check that $\mathrm{h}_{\rm at}^1 (\mathbb{R}^n,\mu) \subset \mathrm{H}_{\rm atb}^1 (\mathbb{R}^n,\mu)$. It is not clear for us what the diagonal part $\mathrm{h}_{\mathrm{diag}}^1(\R^n,\mu)$ should be.
\begin{problem}
Find a diagonal Hardy space $$\mathrm{h}_{\rm diag}^1 (\mathbb{R}^n,\mu)$$ so that the following Davis type decomposition holds
$$\mathrm{H}_{\rm atb}^1 (\mathbb{R}^n,\mu) \, = \, \mathrm{h}_{\rm at}^1 (\mathbb{R}^n,\mu) + \mathrm{h}_{\rm diag}^1 (\mathbb{R}^n,\mu).$$
\end{problem}

\vskip5pt

\noindent \textbf{Acknowledgement.} J.M. Conde-Alonso and J. Parcet are partially supported by the European Research Council ERC StG-256997-CZOSQP, the Spanish grant MTM2010-16518 and by ICMAT Severo Ochoa Grant SEV-2011-0087 (Spain). 

\bibliographystyle{amsplain}

\vskip20pt

\hfill \noindent \textbf{Jose M. Conde-Alonso} \\
\null \hfill Instituto de Ciencias Matem{\'a}ticas \\ \null \hfill
CSIC-UAM-UC3M-UCM \\ \null \hfill Consejo Superior de
Investigaciones Cient{\'\i}ficas \\ \null \hfill C/ Nicol\'as Cabrera 13-15.
28049, Madrid. Spain \\ \null \hfill\texttt{jose.conde@icmat.es}

\vskip5pt


\hfill \noindent \textbf{Javier Parcet} \\
\null \hfill Instituto de Ciencias Matem{\'a}ticas \\ \null \hfill
CSIC-UAM-UC3M-UCM \\ \null \hfill Consejo Superior de
Investigaciones Cient{\'\i}ficas \\ \null \hfill C/ Nicol\'as Cabrera 13-15.
28049, Madrid. Spain \\ \null \hfill\texttt{javier.parcet@icmat.es}
\end{document}